\documentclass[12pt]{amsart}

\usepackage{amssymb}
\usepackage{amscd}

\newcommand{\bt}{\begin{theorem}}
\newcommand{\et}{\end{theorem}}
\newcommand{\bp}{\begin{proposition}}
\newcommand{\ep}{\end{proposition}}
\newcommand{\bq}{\begin{question}}
\newcommand{\eq}{\end{question}}
\newcommand{\bl}{\begin{lemma}}
\newcommand{\el}{\end{lemma}}
\newcommand{\br}{\begin{result}}
\newcommand{\er}{\end{result}}
\newcommand{\be}{\begin{equation}}
\newcommand{\ee}{\end{equation}}
\newcommand{\bc}{\begin{corollary}}
\newcommand{\ec}{\end{corollary}}
\newcommand{\bex}{\begin{example}}
\newcommand{\eex}{\end{example}}

\newtheorem{theorem}{Theorem}[section]
\newtheorem{corollary}[theorem]{Corollary}
\newtheorem{lemma}[theorem]{Lemma}
\newtheorem{proposition}[theorem]{Proposition}
\newtheorem{result}[theorem]{Result}
\newtheorem{example}[theorem]{Example}
\newtheorem{question}[theorem]{Question}

\newcommand{\cPA}{\mathcal{P}\!\mathcal{A}}
\newcommand{\cPSA}{\mathcal{P}\!\mathcal{S}\!\mathcal{A}}
\newcommand{\cAut}{\mathcal{A}\mathit{ut}}
\newcommand{\cB}{\mathcal{B}}

\newcommand{\cD}{\mathcal{D}}
\newcommand{\cH}{\mathcal{H}}
\newcommand{\cI}{\mathcal{I}}
\newcommand{\cJ}{\mathcal{J}}
\newcommand{\cK}{\mathcal{K}}
\newcommand{\cL}{\mathcal{L}}
\newcommand{\cR}{\mathcal{R}}

\newcommand{\cU}{\mathcal{U}}

\newcommand{\tsq}{\,${\tiny$\square$}$\,}
\newcommand{\epr}{\hspace{\fill}$\Box$}

\newcommand{\Ng}{N_}

\newcommand{\med}{\medskip}
\newcommand{\sm}{\smallskip}

\begin{document}

\noindent {\em J. Aust. Math. Soc.} {\bf 81} (2006), 185 -- 198
\vspace{0.05in}\\
{\em arXiv version}: layout, fonts, pagination and numbering of sections and theorems may vary from the JAMS published version
\vspace{0.1in}\\

\title{Inverse semigroups determined by their partial automorphism monoids}
\author{Simon M. Goberstein}
\address{\hspace{-0.205in}Department of Mathematics \!\&\! Statistics,
California State University, Chico, CA 95929, USA, e-mail: sgoberstein@csuchico.edu}
\begin{abstract}
The partial automorphism monoid of an inverse semigroup is an inverse monoid consisting of all isomorphisms between its inverse subsemigroups. We prove that a tightly connected fundamental inverse semigroup $S$ with no isolated nontrivial subgroups is lattice determined ``modulo semilattices'' and if $T$ is an inverse semigroup whose partial automorphism monoid is isomorphic to that of $S$, then either $S$ and $T$ are isomorphic or they are dually isomorphic chains relative to the natural partial order; a similar result holds if $T$ is any semigroup and the inverse monoids consisting of all isomorphisms between subsemigroups of $S$ and $T$, respectively, are isomorphic. Moreover, for these results to hold, the conditions that $S$ be tightly connected and have no isolated nontrivial subgroups are essential. 
\vspace{0.08in}
\end{abstract}
\maketitle
\font\caps=cmcsc10 scaled \magstep1   % Theorems, etc.
\setcounter{section}{-1}

\section{Introduction}
%\vspace{0.05in} 
\sm A {\em partial automorphism} of an algebraic structure $A$ of a certain type is any isomorphism between its substructures (including, if necessary, the empty one), and the set of all partial automorphisms of $A$ with respect to composition is an inverse monoid called the {\em partial automorphism monoid} of $A$. The problem of characterizing algebras of various types by their partial automorphism monoids was posed by Preston in \cite{key18}. Since idempotent partial automorphisms correspond to subalgebras, this is closely related to the problem of characterizing algebras by their subalgebra lattices. In \cite{key4}, we described large classes of combinatorial inverse semigroups determined by their lattices of inverse subsemigroups and partial automorphism monoids (Theorems 5 and 8, respectively) and showed that these theorems do not hold for fundamental inverse semigroups containing isolated nontrivial subgroups (Proposition 10). However, the problem of whether the principal results of \cite{key4} can be extended from combinatorial to fundamental inverse semigroups with {\em no isolated} nontrivial subgroups has remained open. The purpose of this article is to solve this problem in the affirmative. In Section 2, we study lattice isomorphisms of inverse semigroups and show that so-called tightly connected fundamental inverse semigroups without nontrivial isolated subgroups are lattice determined ``modulo semilattices'' (Theorem \ref{24}). Using this theorem, we prove in Section 3 that any tightly connected fundamental inverse semigroup $S$ with no nontrivial isolated subgroups is determined (up to a dual isomorphism if $S$ is a chain with respect to the natural partial order) in the class of {\em all inverse} semigroups by its monoid of isomorphisms between {\em inverse} subsemigroups (Theorem \ref{32}), and in the class of {\em all} semigroups by its  monoid of isomorphisms between {\em arbitrary} subsemigroups (Theorem \ref{34}). Examples described in Propositions 9 and 10 of \cite{key4} show that the conditions imposed on $S$ in Theorems \ref{24}, \ref{32} and \ref{34}, are essential. A few concluding remarks and open questions are contained in Section 4. We use \cite{key1} and \cite{key7} as standard references for the algebraic theory of semigroups, in particular with regard to Green's relations $\cH$, $\cL$, $\cR$, $\cD$ and $\cJ$, and refer to \cite{key16} for an extensive treatment of the theory of inverse semigroups.

The main results of the paper were reported at the Special Session on Semigroup Theory of the $999$th Meeting of the American Mathematical Society held at Vanderbilt University on October 16-17, 2004. 
\sm
\def\bfseries{\normalsize\caps}
\section{Preliminaries}\sm
%\vspace{0.05in}
Denote by $\cI_X$ the symmetric inverse monoid on a set $X$. Let $\varphi$ be a bijection of $X$ onto a set $Y$. For any $\alpha\in\cI_X$, define $\alpha(\varphi\tsq\varphi)=\varphi^{-1} \circ\alpha\circ\varphi$. Clearly $\varphi\tsq\varphi$ is a bijection of $\cI_X$ onto $\cI_Y$, and if $\,\cU$ is a subsemigroup of $\,\cI_X$, then $(\varphi\tsq\varphi)\vert_{\cU}$ is an isomorphism of $\,\cU$ onto $\,\cU(\varphi\tsq\varphi)$.

Let $S$ be an arbitrary semigroup. Denote by $E_S$ the set of idempotents of $S$. For any $x\in S$ and $\cK\in\{\cH,\cL,\cR,\cD,\cJ\}$, denote by $K_x$ the $\cK$-class of $S$ containing $x$, and by 
$J(x)$ the principal two-sided ideal of $S$ generated by $x$. Let $J_x\leq J_y$ if and only if $J(x)\subseteq J(y)$ for $x, y\in S$. Then $\leq$ is a partial order on the set of $\cJ$-classes of $S$. Similarly one can partially order the set of $\cL$-classes and the set of $\cR$-classes of $S$. If $U$ is a subsemigroup of $S$, to distinguish its Green's relations from those on $S$, we will use superscripts. If $U$ is a regular (in particular, inverse) subsemigroup of $S$, then $\cK^U=\cK^S\cap(U\times U)$ for $\cK\in\{\cH,\cL,\cR\}$ (see Hall \cite[Result 9]{key6}).

Let $S$ be an inverse semigroup. We say that $x\in S$ is a {\em group element} if it belongs to some subgroup of $S$; otherwise $x$ is a {\em nongroup element}. Denote by $\Ng S$ the set of all nongroup elements of $S$. Recall that $S$ is termed {\em combinatorial} \cite{key16} if $\cH=1_S$, that is, if every nonidempotent element of $S$ is nongroup. If $x\in S$ and $H_x=\{x\}$, the $\cD$-class $D_x$ will be called {\em combinatorial}. Following Jones \cite{key9}, we say that an idempotent $e$ of $S$ (and each subgroup of $H_e$) is {\em isolated} if $D_e=H_e$, and {\em nonisolated} otherwise. For any $X\subseteq S$, denote by $\langle X \rangle$ the inverse subsemigroup of $S$ generated by $X$. If 
$x\in S$, we say that $\langle x\rangle$ is a {\em monogenic inverse subsemigroup} of $S$ generated by $x$, and if $S=\langle x\rangle$, the inverse semigroup $S$ is {\em monogenic}. A monogenic inverse semigroup $\langle x\rangle$ such that $xx^{-1} > x^{-1}x$ is an inverse monoid with identity $xx^{-1}$ called the {\em bicyclic semigroup}; we denote it by $\cB(x,x^{-1})$. It is well known (see \cite[Theorem 2.53]{key1}) that $\cB(x,x^{-1})$ consists of a single $\cD$-class and its idempotents form a chain: $1=xx^{-1}>x^{-1}x>x^{-2}x^2>\cdots$. If $S$ contains no bicyclic subsemigroup, it is called {\em completely semisimple}. The structure of monogenic inverse semigroups is described in detail in \cite[Chapter IX]{key16}. We recall only a few basic facts about them. 

Let $S=\langle x\rangle$ be a monogenic inverse semigroup. Then $\cD=\cJ$ and the partially ordered set of $\cD$-classes (=$\cJ$-classes) of $S$ is a chain with the largest element $D_x$. It is obvious that one of the following holds: (a) $xx^{-1}=x^{-1}x$, (b) $xx^{-1}$ and $x^{-1}x$ are incomparable with respect to the natural partial order, (c) $xx^{-1}>x^{-1}x$ or $x^{-1}x>xx^{-1}$. In case (a), $S=D_x$ is a cyclic group. In case (b), $D_x=\{x,x^{-1},xx^{-1},x^{-1}x\}$ is the greatest $\cD$-class of $S$ (and, of course, $S\setminus D_x$ is an ideal of $S$). Finally, in case (c), $S=D_x$ is  bicyclic, it is either $\cB(x,x^{-1})$ or $\cB(x^{-1},x)$. In case (b), either $S$ is a free monogenic inverse semigroup or it contains a smallest ideal $K$, the so-called {\em kernel} of $S$; in the latter case, either $K$ is a bicyclic semigroup or a cyclic group. Note that $D_a$ is combinatorial for any $a\in S$ except for the case when $a\in K$ and $K$ is a group. 

An inverse semigroup $S$ is called {\em fundamental} if $1_S$ is the only idempotent-separating congruence on $S$. Since all idempotent-separating congruences on a regular semigroup are contained in $\cH$, it is immediate that combinatorial inverse semigroups are fundamental (but clearly not conversely). Fundamental inverse semigroups, introduced by Munn \cite{key14} (and independently by Wagner \cite{key24}), constitute one of the most important classes of inverse semigroups. Munn, in particular, provided an effective method of describing all fundamental inverse semigroups with the given semilattice of idempotents. Let $E$ be a semilattice. Then the {\em Munn semigroup} $T_E$ of $E$ is an inverse semigroup (under composition) consisting of all isomorphisms between principal ideals of $E$ (see \cite[Section V.4]{key7}). It is common to identify each $e\in E$ with $1_{Ee}\in T_E$, so the semilattice of idempotents of $T_E$ is identified with $E$. If $S$ is an inverse semigroup, a subset $K$ of $S$ is called {\em full} if $E_S\subseteq K$. Munn proved (see \cite[Theorem 2.6]{key14}) that an inverse semigroup $S$ with $E_S=E$ is fundamental if and only if $S$ is isomorphic to a full inverse subsemigroup of $T_E$; in particular, $T_E$ itself is fundamental. It is well known that $S$ is fundamental if and only if for any $x, y\in S$, if $x^{-1}ex=y^{-1}ey$ for all $e\in E_S$, then $x=y$ (see \cite[Section V.3]{key7}). This result was used in the proof of Lemma 2.1 of \cite{key3}, which provides a convenient criterion for a bijection between an arbitrary inverse semigroup and a fundamental one to be an isomorphism.
\br\label{11}{\rm (From \cite[Lemma 2.1]{key3})} Let $S$ and $T$ be inverse semigroups and $\varphi$ a bijection of $S$ onto $T$. If $S$ is fundamental and $\varphi\vert_{E_S}$ is an isomorphism of $E_S$ onto $E_T$, then $\varphi$ is an isomorphism of $S$ onto $T$ if and only if $(s^{-1}es)\varphi=(s\varphi)^{-1}(e\varphi)(s\varphi)$ for all $s\in S$ and $e\in E_S$.
\er
By modifying the proof of Lemma 2.1 of \cite{key3}, we also obtain the following 
\bl\label{12} Let $S$ and $T$ be inverse semigroups and $\varphi$ a bijection of $S$ onto $T$, preserving $\cL$-classes, such that $\varphi\vert_{E_S}$ is an isomorphism of $E_S$ onto $E_T$. Suppose that $(fx)\varphi=(f\varphi)(x\varphi)$ for all $x\in\Ng S$ and $f\leq xx^{-1}$. Then $(x^{-1}ex)\varphi=(x\varphi)^{-1} (e\varphi)(x\varphi)$ for all $e\in E_S$ and $x\in \Ng S$.
\el 
{\bf Proof.} Let $x\in\Ng S$ and $e\in E_S$. By assumption, $(fx)\varphi=(f\varphi)(x\varphi)$ whenever $f\leq xx^{-1}$, so
\((ex)\varphi=(exx^{-1})\varphi\cdot x\varphi=e\varphi\cdot (xx^{-1})\varphi\cdot x\varphi=
(e\varphi)(x\varphi).\) Since $\varphi$ preserves $\cL$-classes, $(s^{-1}s)\varphi=(s\varphi)^{-1}(s\varphi)$ for each $s\in S$ because  
$s^{-1}s\cL s$ and $(s\varphi)^{-1}(s\varphi)$ is the only idempotent in $L_{s\varphi}$. Thus
\((x^{-1}ex)\varphi=[(ex)\varphi]^{-1}(ex)\varphi=(e\varphi\cdot x\varphi)^{-1}(e\varphi\cdot x\varphi)=(x\varphi)^{-1}(e\varphi)(x\varphi).\)\epr\\

The following auxiliary result is established by applying an argument from the second paragraph of the proof of Lemma 2 of \cite{key4} to a slightly more general situation. 
\br\label{13}{\rm (From \cite[the proof of Lemma 2]{key4})} Let $S$ be an inverse semigroup, $x\in\Ng S$, and 
$e\in E_S$. If $u\in\langle e, x\rangle\cap R_e$ and $u\not=e$, then $u=ex^m$ for some nonzero integer $m$.
\er
\sm
\section{Lattice determinability}
\sm Let $S$ be an inverse semigroup. To indicate that $H$ is an inverse subsemigroup of $S$, we write $H\leq S$. Since 
$\emptyset\leq S$, the set of all inverse subsemigroups of $S$, partially ordered by inclusion, is a complete (and compactly generated) lattice denoted by $\cL(S)$ (as in \cite{key9}). Let $T$ be an inverse semigroup such that there is an isomorphism $\Psi$ of $\cL(S)$ onto $\cL(T)$. Then $S$ and $T$ are called {\em lattice isomorphic}, and $\Psi$ is a {\em lattice isomorphism} of $S$ onto $T$. We say that a mapping $\psi : S\rightarrow T$ {\em induces} $\Psi$ (or $\Psi$ is {\em induced} by $\psi$) if $H\Psi=H\psi$ for all $H\leq S$. If $S$ is isomorphic to every inverse semigroup that is lattice isomorphic to $S$, then $S$ is called {\em lattice determined}, and if each lattice isomorphism of $S$ onto an inverse semigroup $T$ is induced by an isomorphism of $S$ onto $T$, then $S$ is {\em strongly lattice determined}.

Let $S$ and $T$ be inverse semigroups and $\Psi$ a lattice isomorphism of $S$ onto $T$. It is clear that an inverse subsemigroup $U$ of $S$ is an atom of $\cL(S)$ if and only if $U=\{e\}$ for some $e\in E_S$. Thus there is a unique bijection $\psi_E$ of $E_S$ onto $E_T$ defined by the formula $\{e\}\Psi=\{e\psi_E\}$ for all $e\in E_S$, and we will say that $\psi_E$ is the {\em $E$-bijection} associated with $\Psi$. Recall that if $X$ is a partially ordered set and 
$x, y\in X$, then $x\parallel y$ means that $x$ and $y$ are incomparable in $X$, and $x\nparallel y$ denotes the negation of $x\parallel y$. It is well known (see \cite[Subsection 36.6 and Chapter XIV, Introduction]{key21}) that for all $e, f\in E_S$, we have $e\nparallel f$ if and only if $e\psi_E\nparallel f\psi_E$, and if $e\parallel f$, then $(ef)\psi_E=(e\psi_E)(f\psi_E)$, which is expressed by saying that $\psi_E$ is a {\em weak isomorphism} of $E_S$ onto $E_T$. 
\br\label{20}{\rm (See \cite[Proposition 1.6 and Corollary 1.7]{key9})} If $S$ and $T$ are inverse semigroups and $\Psi$ a lattice isomorphism of $S$ onto $T$, there is a (unique) bijection $\psi : \Ng S\cup E_S\rightarrow \Ng T\cup E_T$ with the following properties{\rm :}\\
{\rm (a)} $\psi$ extends $\psi_E$, that is, $\psi\vert_{E_S}=\psi_E${\rm ;}\\
{\rm (b)} $\psi$ and $\psi^{-1}$ preserve $\cR$- and $\cL$-classes{\rm ;}\\
{\rm (c)} for every $x\in \Ng S\cup E_S$, we have $\langle x\rangle\Psi=\langle x\psi\rangle$ so, in particular, $(x^{-1})\psi=(x\psi)^{-1}${\rm ;}\\{\rm (d)} if a homomorphism $\gamma : S\rightarrow T$ induces $\Psi$, then $x\psi=x\gamma$ for all $x\in \Ng S\cup E_S$. 
\er

Using the terminology of \cite{key21}, we say that the bijection $\psi : \Ng S\cup E_S\rightarrow \Ng T\cup E_T$ in Result
\ref{20} is the {\em base partial bijection} associated with the lattice isomorphism $\Psi$ of $S$ onto $T$. In the notation of Result \ref{20}, suppose that $S$ has no nontrivial isolated subgroups. Then $T$ also has no nontrivial isolated subgroups (see \cite[Corollary 1.9]{key9}). As shown by Ershova \cite{key2}, in this case there is a bijection $\hat\psi$ of $S$ onto $T$ which extends $\psi$ and retains a number of its properties. This bijection $\hat\psi$ can be constructed as follows. First, set $x\hat\psi=x\psi$ if $x$ is a nongroup element or an isolated idempotent of $S$. Now for every nonisolated idempotent $e$ of $S$, choose and fix an element $r_e\in \Ng S\cap R_e$. It is easily seen that for each $a\in H_e$, there is a unique $q\in H_{r_e}$ such that $a=r_eq^{-1}$,  and we put $a\hat\psi=(r_e\psi)(q\psi)^{-1}$. Then it can be shown that $e\hat\psi=e\psi$ for every $e\in E_S$, $\hat\psi$ preserves $\cL$- and $\cR$-classes, that is, $(ss^{-1})\hat\psi=(s\hat\psi)(s\hat\psi)^{-1}$ and $(s^{-1}s)\hat\psi=(s\hat\psi)^{-1}(s\hat\psi)$ for all $s\in S$, and if $\Psi$ is induced by an isomorphism $\gamma$ of $S$ onto $T$, then $\gamma=\hat\psi$ (see  \cite[Lemmas 1 and 2]{key2}, or \cite[Subsection 43.7]{key21}). Again using the terminology of \cite{key21}, we call $\hat\psi$ the {\em base bijection} of $S$ onto $T$ associated with $\Psi$ (clearly, it depends also on the choice of $r_e\in N_S\cap R_e$ for each nonisolated $e\in E_S$). Of course, if $S$ is combinatorial, then $\hat\psi=\psi$, that is, $\psi$ is the base bijection of $S$ onto $T$.

Let $S$ be an inverse semigroup. If $x\in S$ and $e\in E_S$ are such that $e<xx^{-1}$ and no 
$f\in E_{\langle x\rangle}$ satisfies  $e<f<xx^{-1}$, we say that $e$ is $x${\em -covered} by $xx^{-1}$ and write $e\prec_x xx^{-1}$. Take $x\in S$ and $e\in E_S$ with $e<xx^{-1}$. If for some positive integer $n$, there exist $e_{0}, e_{1}, \ldots, e_n\in E_S$ such that $e=e_{0}<e_{1}<\cdots<e_n=xx^{-1}$ and for every $k=1,\ldots,n$, the idempotent $e_{k-1}$ is $x_{k}$-covered by $e_{k}$ where $x_{k}=e_{k}x$ (so that $x_{k} x^{-1}_{k}=e_{k}$), then 
$(e_{0},e_{1},\ldots,e_n)$ is called a {\em short bypass} from $e$ to $xx^{-1}$ (it is plain that if $x_k\in E_S$, then 
$x_i\in E_S$ for all $i < k$). If for all $x\in S$ and $e\in E_S$ such that $e<xx^{-1}$ there is a short bypass from $e$ to $xx^{-1}$, then $S$ is called {\em shortly connected}. If $x\in\Ng S\cup E_S$ and $e\prec_x xx^{-1}$, we say that $e$ is {\em tightly $x$-covered} by $xx^{-1}$ if either $ex\in\Ng S$ or $ex\in E_S$ (in the latter case, of course, $ex=e$ because $ex\in R_e$). Let $x\in\Ng S\cup E_S$ and $e <xx^{-1}$. If there is a short bypass $(e_{0},e_{1},\ldots,e_n)$ from $e\,(=e_0)$ to $xx^{-1}\,(=e_n)$ such that for every $k=1,\ldots,n$, the idempotent $e_{k-1}$ is tightly $x_{k}$-covered by $e_{k}$ (where, as above, $x_{k}=e_{k}x$), then $(e_{0},e_{1},\ldots,e_n)$ will be called a {\em tight bypass} from $e$ to $xx^{-1}$ (in which case, if $x\in\Ng S$, there is a smallest $m\in\{1, \ldots, n\}$ satisfying $x_m\in\Ng S$, so $x_i\in E_S$ for all $i < m$, and $x_j\in\Ng S$ for all $m\leq j\leq n$). We say that $S$ is {\em tightly connected} if for all $x\in\Ng S\cup E_S$ and $e<xx^{-1}$, there is a tight bypass from $e$ to $xx^{-1}$.

It is obvious that every tightly connected inverse semigroup is shortly connected, and a combinatorial inverse semigroup is shortly connected if and only if it is tightly connected. At the same time, tightly connected fundamental inverse semigroups need not be combinatorial -- see, for instance, Example 2 of \cite{key4}. Although the semigroups in that example have nontrivial isolated subgroups, there are many noncombinatorial tightly connected fundamental inverse semigroups with no nontrivial isolated subgroups. The smallest such example can be constructed as follows. Let $E=\{e_0, e_1, f_0, f_1, f_2, 0\}$ be the semilattice given by the diagram in Figure 1.
\vspace{0.1in}\\
\begin{picture}(300,135)
\put(260,70){\line(-6,-5){60}} \put(140,70){\line(6,-5){60}} \put(170,120){\line(3,-5){30}} \put(170,120){\line(-3,-5){30}}
\put(230,120){\line(3,-5){30}} \put(230,120){\line(-3,-5){30}}  \put(200,70){\line(0,-1){50}}
\put(230,120){\circle*{4}}\put(170,120){\circle*{4}} \put(140,70){\circle*{4}} \put(260,70){\circle*{4}}
\put(200,70){\circle*{4}} \put(200,20){\circle*{4}} \put(203,10){$0$} \put(155,120){$e_0$} \put(235,120){$e_1$}
\put(125,68){$f_0$} \put(205,68){$f_1$} \put(265,68){$f_2$}\put(180,-5){$\text{Figure }1$}
\end{picture}
\vspace{0.2in}\\
Let $S=T_E$ be the Munn semigroup of $E$. It is easily seen that $S$ is tightly connected and consists of three $\cD$-classes: $D_0=\{0\}$, $D_{f_0}\;(=D_{f_1}=D_{f_2})$, and $D_{e_0}\;(=D_{e_1})$, such that $D_0<D_{f_0}<D_{e_0}$. Moreover, $\{0\}\cup D_{f_0}$ is a $10$-element combinatorial inverse semigroup, and $D_{e_0}$, the top $\cD$-class of $S$, consists of four nontrivial $\cH$-classes. Thus $S$ is a noncombinatorial tightly connected fundamental inverse semigroup with no nontrivial isolated subgroups.
\bl\label{21} Let $S$ be an inverse semigroup, and let $U=\langle e, x\rangle$ for some $x\in\Ng S$ and $e\in E_S$ such that $ex\in\Ng S$ and $e$ is $x$-covered by $xx^{-1}$. Then $H^U_e=\{e\}$.
\el
{\bf Proof.} Suppose $u\in H^U_e$ and $u\not=e$. By Result \ref{13}, $u=ex^m$ for some $m\not=0$. Then $ex^mx^{-m}=e$ and $x^{-m}ex^m=e$, so $x^me=ex^m$, $x^{-m}e=ex^{-m}$, and $ex^{-m}x^m=e$. Let $n=|m|$. Clearly $m\not=-1$ (otherwise $ex=(x^{-1}e)^{-1}=(ex^{-1})^{-1}\in H^U_e$, a contradiction). Thus $n\geq 2$. If $xx^{-1}>x^{-1}x$, then $e=ex^{-n}x^n\leq x^{-2}x^2<x^{-1}x<xx^{-1}$, contradicting $e\prec_x xx^{-1}$. Hence $xx^{-1}\not>x^{-1}x$, in which case $e=ex^nx^{-n}\leq x^2x^{-2} < xx^{-1}$, so we have $e=x^2x^{-2}\in\langle x\rangle$. Thus $U=\langle x\rangle$ and $D^U_e$ is a combinatorial $\cD$-class of $U$ since $ex\in D^U_e$ and $ex\in\Ng S$. This contradicts the assumption that $e\not=u\in H^U_e$. Therefore $H^U_e=\{e\}$.\epr\\ 
 
For combinatorial inverse semigroups the following theorem was proved in \cite{key4}: 
\br\label{22} {\rm (\cite[Theorem 5]{key4})} Let $S$ be a combinatorial inverse semigroup, $T$ an inverse semigroup and $\Psi$ a lattice isomorphism of 
$S$ onto $T$. Let $\psi$ be the base bijection of $S$ onto $T$ associated with $\Psi$ (so, in particular,
$\psi_{E}=\psi\vert_{E_S}$). Suppose that $S$ is shortly connected (equivalently, tightly connected) and $\psi_{E}$ is an
isomorphism of $E_S$ onto $E_T$. Then $\psi$ is the unique isomorphism of $S$ onto $T$ which induces $\Psi$. \er

We are going to extend this theorem to the class of tightly connected fundamental inverse semigroups with no nontrivial isolated subgroups. A key role in the proof of Result \ref{22} was played by Lemma 2 of \cite{key4}, which can be modified to establish the following more general result: 
\bl\label{23} Let $S$ be an inverse semigroup, let $\Psi$ be a lattice isomorphism of $S$ onto an inverse semigroup $T$, and let $\psi : N_S\cup E_S\rightarrow N_T\cup E_T$ be the base partial bijection associated with $\Psi$. Suppose that 
$\psi\vert_{E_S}$ is an isomorphism of $E_S$ onto $E_T$. Then $(ex)\psi=(e\psi)(x\psi)$ for any $x\in \Ng S\cup E_S$ and any $e\in E_S$ such that $e$ is tightly $x$-covered by $xx^{-1}$.
\el 
{\bf Proof.} In the following proof, reference to \cite{key4} means with respect to parts of the proof (almost verbatim) of Lemma 2 of \cite{key4}. Let $x\in\Ng S\cup E_S$ and $e\in E_S$ be such that $e$ is tightly $x$-covered by $xx^{-1}$. Denote $U=\langle e, x\rangle$ and $V=\langle e\psi, x\psi\rangle$. Then the restriction of $\Psi$ to $\cL(U)$ is a lattice isomorphism of $U$ onto $U\Psi=\langle e\psi\rangle\vee\langle x\psi\rangle=V$. We show that $(ex)\psi=(e\psi)(x\psi)$. This holds for $x\in E_S$, so assume that $x\in \Ng S$. By 
\cite{key4}, we may suppose that $xx^{-1}\nless x^{-1}x$. By \cite{key4} also, $(ex)\psi\cR(e\psi)(x\psi)$. Since $e$ is tightly $x$-covered by $xx^{-1}$, either $ex=e$ or $ex\in\Ng S$.

{\bf Case I}. Suppose $ex=e$. Then it is easily seen that $e$ is a zero for all elements of $\langle x\rangle$ so, in particular, $e\leq x^{-1}x^2x^{-1}$. Note that $x^{-1}x^2x^{-1} < xx^{-1}$ since $xx^{-1}\not<x^{-1}x$. It follows that $e=x^{-1}x^2x^{-1}$. Then $e=x^2x^{-1}ex^{-1}x^2=x^2$, so $U=\langle x\rangle=\{x, x^{-1}, xx^{-1}, x^{-1}x, e\}$ is a five-element Brandt semigroup \cite{key16} which is strongly lattice determined (this can be easily shown directly and also follows from  \cite[Theorem 42.4]{key21}). Therefore $(ex)\psi=e\psi=(e\psi)(x\psi)$. 

{\bf Case II}. Suppose $ex\in\Ng S$. By Lemma \ref{21}, $H^U_e=\{e\}$ and so $H^V_{e\psi}=\{e\psi\}$ since $\Psi\vert_{\cL(H^U_e)}$ is a lattice isomorphism of $H^U_e$ onto $H^V_{e\psi}$ (see \cite[Corollary 1.2]{key8}). Hence $(e\psi)(x\psi)\in\Ng T$ or $(e\psi)(x\psi)=e\psi$. If $(e\psi)(x\psi)=e\psi$, as in Case I, $e\psi$ is a zero for all elements of $\langle x\psi\rangle$, which implies $e\leq f$ for all $f\in E_{\langle x\rangle}$, and $e\in\langle x\rangle=U$, whence $ex\in D^U_e=H^U_e$, a contradiction. Thus $(e\psi)(x\psi)\in\Ng T$, so that $(e\psi)(x\psi)=u\psi$ where $u\in R^U_e$ and $u\not=e$. By Result \ref{13}, $u=ex^n$ or $=ex^{-n}$ for some integer $n\geq 1$. If $u=ex$, we are done. In all other cases, as in \cite{key4}, we have $e=x^{-1}x^2x^{-1}\in\langle x\rangle$. If $\langle x\rangle$ is combinatorial, $(ex)\psi=(e\psi)(x\psi)$ because combinatorial monogenic inverse semigroups are strongly lattice determined (see \cite[Theorems 42.2 and 42.4]{key21}). Otherwise the kernel of $\langle x\rangle$ is a nontrivial (cyclic) group which does not contain $ex$, so $(ex)\psi=(e\psi)(x\psi)$ by \cite[Lemma 42.7]{key21}.\epr\\

We are ready to prove the main result of this section. 
\bt\label{24} Let $S$ be a tightly connected fundamental inverse semigroup with no nontrivial isolated subgroups, $T$ an inverse semigroup, and $\Psi$ a lattice isomorphism of $S$ onto $T$. As above, let $\psi_E$ denote the $E$-bijection, $\psi$ the base partial bijection, and $\hat\psi$ the base bijection of $S$ onto $T$ associated with $\Psi$. Suppose that $\psi_E$ is an isomorphism of $E_S$ onto $E_T$. Then $\hat\psi$ is an isomorphism of $S$ onto $T$. If $S$ is combinatorial or completely semisimple, then $\hat\psi$ is the unique isomorphism of $S$ onto $T$ inducing $\Psi$.
\et 
{\bf Proof.} Since $\hat\psi\vert_{\Ng S\cup E_S}=\psi$, we will write $s\psi$ instead of $s\hat\psi$ for any $s\in\Ng S\cup E_S$. Take arbitrary $x\in\Ng S$ and $f < xx^{-1}$. By assumption, there is a tight bypass $(f_{0},f_{1},\ldots, f_n)$ from $f=f_0$ to $xx^{-1}=f_n$ for some $n\geq 1$. For $0\leq i\leq n$, set $x_i=f_i x$. Then $x_ix^{-1}_i=f_i$ for all $0\leq i\leq n$. Moreover, $x_{k-1}=f_{k-1}x_k$ and $f_{k-1}$ is tightly $x_k$-covered by $f_k$ for each $1\leq k\leq n$. Since $x_n=x\in\Ng S$, there is a smallest $m\in\{1, \ldots, n\}$ satisfying $x_m\in\Ng S$, so $x_i\in E_S$ for all $i < m$, and $x_j\in\Ng S$ for all $m\leq j\leq n$. Since $\psi_E$ is an isomorphism of $E_S$ onto $E_T$, \((fx)\psi=(f_0f_{m-1}x)\psi=(f_0x_{m-1})\psi=(f_0\psi)(x_{m-1}\psi)\).  Now using Lemma \ref{23}, we obtain
\begin{eqnarray*}
x_{m-1}\psi &=& (f_{m-1}x_m)\psi=(f_{m-1}\psi)(x_m\psi)=(f_{m-1}\psi)\cdot(f_mx_{m+1})\psi=\cdots\\
                    &=& (f_{m-1}\psi)\cdot(f_{n-1}x_n)\psi=(f_{m-1}\psi)(f_{n-1}\psi)(x\psi)=(f_{m-1}\psi)(x\psi).
\end{eqnarray*}
It follows that $(fx)\psi=(f_0\psi)(f_{m-1}\psi)(x\psi)=(f\psi)(x\psi)$, and since $\psi$ preserves $\cR$-classes, we also have
$(xx^{-1}\cdot x)\psi=x\psi=(x\psi)(x\psi)^{-1}(x\psi)=(xx^{-1})\psi\cdot x\psi$. Therefore $(gx)\psi=(g\psi)(x\psi)$ for all $g\leq xx^{-1}$, and hence, by Lemma \ref{12}, $(x^{-1}ex)\psi=(x\psi)^{-1} (e\psi)(x\psi)$ for all $e\in E_S$.

Now take an arbitrary group element $a\in S$. Let us show that $(a^{-1}ea)\hat\psi=(a\hat\psi)^{-1}(e\hat\psi)(a\hat\psi)$ for all $e\in E_S$. If $a\in E_S$, this holds because $\psi_E$ is an isomorphism of $E_S$ onto $E_T$. Thus we may assume that $a\not\in E_S$. Denote $aa^{-1}$ by $f$. Recall that we have fixed an element $r\!_{_f}\in \Ng S \cap R_f$ and that there is a unique $q\in H_{r\!_{_f}}$ such that $a=r\!_{_f} q^{-1}$. To shorten notation, set $r=r\!_{_f}$. Take an arbitrary $e\in E_S$. It follows from the first paragraph of the proof that $(q\cdot r^{-1}er\cdot q^{-1})\psi=(q\psi)\cdot (r^{-1}er)\psi\cdot (q\psi)^{-1}$ and $(r^{-1}er)\psi=(r\psi)^{-1}(e\psi)(r\psi)$ because $r$ and $q$ are nongroup elements of $S$ and $(q^{-1})\psi=(q\psi)^{-1}$. Therefore
\begin{eqnarray*}
(a^{-1}ea)\hat\psi &=& ((rq^{-1})^{-1}\cdot e\cdot rq^{-1})\psi=(q\cdot r^{-1} e r\cdot
q^{-1})\psi\\
&=&q\psi\cdot(r^{-1}er)\psi\cdot(q\psi)^{-1}=q\psi\cdot [(r\psi)^{-1}\cdot e\psi\cdot r\psi]\cdot
(q\psi)^{-1}\\
&=&[r\psi\cdot(q\psi)^{-1}]^{-1}\cdot e\psi\cdot[r\psi\cdot(q\psi)^{-1}]=(a\hat\psi)^{-1}(e\hat\psi)(a\hat\psi).
\end{eqnarray*}
We have shown that $(s^{-1}es)\hat\psi=(s\hat\psi)^{-1}(e\hat\psi)(s\hat\psi)$ for all $s\in S$ and all $e\in E_S$. Since $S$ is a fundamental inverse semigroup, according to Result \ref{11}, $\hat\psi$ is an isomorphism of $S$ onto $T$.

If $S$ is combinatorial, then $\hat\psi=\psi$ and, by Result \ref{22}, $\psi$ is the unique isomorphism of $S$ onto $T$ inducing $\Psi$. Suppose that $S$ is completely semisimple. By  \cite[Lemma 2.4]{key10}, $\Psi$ is induced by $\hat\psi$. Hence, in view of Ershova's result cited earlier (see \cite[Proposition 43.7.3]{key21}), $\hat\psi$ is the unique isomorphism of $S$ onto $T$ inducing $\Psi$. This completes the proof.\epr
\sm
\section{$\cPA$-determinability}
\sm Let $S$ be an inverse semigroup. In this paper, we define a {\em partial automorphism} of $S$ to be any isomorphism
between its inverse subsemigroups and denote by $\cPA(S)$ the set of all partial automorphisms of $S$. It is easy to see that with respect to composition $\cPA(S)$ is an inverse submonoid of $\cI_S$. We call $\cPA(S)$ the {\em partial automorphism monoid} of $S$. The group of units of $\cPA(S)$ is $\cAut(S)$, the automorphism group of $S$, and the semilattice of idempotents of $\cPA(S)$ is a lattice isomorphic to $\cL(S)$.

Let $S$ and $T$ be inverse semigroups. If $\cPA(S)\cong\cPA(T)$, then $S$ and $T$ are said to be $\cPA${\em-isomorphic}, and any isomorphism of $\cPA(S)$ onto $\cPA(T)$ is called a $\cPA${\em-isomorphism} of $S$ onto $T$. Let $\Phi$ be a
$\cPA$-isomorphism of $S$ onto $T$. We say that $\Phi$ is {\em induced} by a bijection $\varphi\!:S\rightarrow T$ if $\Phi=(\varphi\tsq\varphi)\vert_{\cPA(S)}$, that is, if for all 
$\alpha\in\cPA(S)$ and $x,y\in S$, we have $\;x\alpha=y$ if and only if $(x\varphi)(\alpha\Phi)=y\varphi$. Let $\xi$ be any bijection of $S$ onto $T$. It is clear that $(\xi\tsq\xi)\vert_{\cPA(S)}$ is a $\cPA$-isomorphism of $S$ onto $T$ precisely when $\cPA(S)(\xi\tsq\xi)=\cPA(T)$. In particular, any isomorphism (or antiisomorphism) of $S$ onto $T$ induces a $\cPA$-isomorphism of $S$ onto $T$. An inverse semigroup $S$ is called $\cPA${\em-determined} if it is isomorphic to any inverse semigroup $\cPA$-isomorphic to $S$, and {\em strongly $\cPA$-determined} if each $\cPA$-isomorphism of $S$ onto an inverse semigroup $T$ is induced by an isomorphism of $S$ onto $T$.

Let $S$ and $T$ be $\cPA$-isomorphic inverse semigroups and $\Phi$ a $\cPA$-isomorphism of $S$ onto $T$. For any 
$H\leq S$, define $H\Phi^*$ by the formula $1_H\Phi=1_{H\Phi^*}$. Then $\Phi^*$ is a lattice isomorphism of $S$ onto $T$. We will denote by $\varphi_E$ the $E$-bijection and by $\varphi$ the base partial bijection associated with $\Phi^*$, and say that $\varphi_E$ and $\varphi$ are {\em associated with $\Phi$}. As shown in Section 2, if $S$ has no nontrivial isolated subgroups, we can extend $\varphi$ to the base bijection $\hat\varphi$ of $S$ onto $T$ associated with $\Phi^*$, and again we will say that $\hat\varphi$ is {\em associated with} the $\cPA$-isomorphism $\Phi$.

\br\label{31}{\rm (A corollary to \cite[Theorem]{key19})} Let $S$ be a semilattice and $T$ an inverse semigroup. Then $\cPA(S)\cong\cPA(T)$ if and only if either $S\cong T$ or $S$ is a chain and $T\cong S^d$. Moreover, any $\cPA$-isomorphism $\Phi$ of $S$ onto $T$ is induced by the $E$-bijection $\varphi_{E}$ associated with $\Phi$, and $\varphi_{E}$ is either an isomorphism or, if $S$ is a chain and $T\cong S^d$, a dual isomorphism of $S$ onto $T$. 
\er 
Using Result \ref{31} and Theorem \ref{24}, we can easily prove the following theorem which establishes $\cPA$-determinability of tightly connected fundamental inverse semigroups having no nontrivial isolated subgroups (with the exception of chains that are $\cPA$-determined up to a dual isomorphism) and thus extends Theorem 8 of \cite{key4}.

\bt\label{32} Let $S$ be a tightly connected fundamental inverse semigroup with no nontrivial isolated subgroups and $T$ an
inverse semigroup. Then $\cPA(S)\cong\cPA(T)$ if and only if either $S\cong T$ or $(S, \leq)$ and $(T, \leq)$ are dually
isomorphic chains. More specifically, let $\Phi$ be a $\cPA$-isomorphism of $S$ onto $T$. As above, denote by $\varphi$ the base partial bijection and by $\hat\varphi$ the base bijection of $S$ onto $T$ associated with $\Phi$. Then either $(S, \leq)$ and $(T, \leq)$ are dually isomorphic chains and $\varphi$ is the unique dual isomorphism of $(S, \leq)$ onto $(T, \leq)$ inducing $\Phi$ or $\hat\varphi$ is an isomorphism of $S$ onto $T$. If $S$ is combinatorial or completely semisimple, then $\hat\varphi$ is the unique isomorphism of $S$ onto $T$ inducing $\Phi$.
\et 
{\bf Proof.} Let $\Phi$ be a $\cPA$-isomorphism of $S$ onto $T$. According to \cite[Lemma 7]{key4}, the restriction of $\Phi$ to $\cPA(E_S)$ is a $\cPA$-isomorphism of $E_S$ onto $E_T$. Hence, by Result \ref{31}, either $\varphi_E$ is an isomorphism of $E_S$ onto $E_T$, or $(E_S, \leq)$ and $(E_T, \leq)$ are dually isomorphic chains and $\varphi_E$ is a dual isomorphism of $E_S$ onto $E_T$. If the latter holds, then according to the argument in the last paragraph of the proof of Theorem 8 of \cite{key4}, we have $S=E_S$ and $T=E_T$, so that $(S, \leq)$ and $(T, \leq)$ are dually isomorphic chains, and $\varphi\,(=\varphi_E)$ is the unique dual isomorphism of $(S, \leq)$ onto $(T, \leq)$ inducing $\Phi$.

Now suppose that $\varphi_E$ is an isomorphism of $E_S$ onto $E_T$. Then, by Theorem \ref{24}, $\hat\varphi$ is an isomorphism of $S$ onto $T$. If $S$ is combinatorial, then $\hat\varphi=\varphi$ and, by \cite[Theorem 8]{key4}, $\varphi$ is the unique isomorphism of $S$ onto $T$ inducing $\Phi$. If $S$ is completely semisimple, according to Theorem \ref{24}, $\hat\varphi$ induces $\Phi^*$. Similarly to the combinatorial case, this implies that $\hat\varphi$ induces $\Phi$ (see the corresponding part of the proof of \cite[Theorem 8]{key4}). This completes the proof.\epr\\

In \cite[Proposition 9]{key4}, we constructed examples of $\cPA$-isomorphic (and thus lattice isomorphic) completely
semisimple {\em combinatorial} inverse semigroups which are not isomorphic, thereby showing that the requirement that $S$ in Theorems 5 and 8 of \cite{key4} be shortly connected is essential. Of course, the same examples show that the requirement that $S$ in Theorems \ref{24} and \ref{32} be tightly connected is essential as well. Furthermore, by Proposition 10 of \cite{key4}, there exist finite (and thus shortly connected) {\em fundamental} inverse semigroups, containing nontrivial isolated subgroups, which are $\cPA$-isomorphic but not isomorphic. It is easily seen that the inverse semigroups, constructed in \cite[Proposition 10]{key4}, are tightly connected. This shows that the requirement that a tightly connected fundamental inverse semigroup $S$ in Theorems \ref{24} and \ref{32} have no nontrivial isolated subgroups is also essential.

We turn now to the problem of determinability of inverse semigroups by partial automorphism monoids in the class of all
semigroups. Let $S$ be any (not necessarily inverse) semigroup. We will say that an isomorphism between subsemigroups of $S$ is a {\em partial $s$-automorphism} of $S$. (We regard $\emptyset$ as a subsemigroup of $S$, so it is also a partial 
$s$-automorphism of $S$.) Let $\cPSA(S)$ denote the set of all partial $s$-automorphisms of $S$. It is easily seen that with respect to composition $\cPSA(S)$ is an inverse submonoid of $\cI_S$ whose semilattice of idempotents is isomorphic to the lattice of all subsemigroups of $S$. If $T$ is a semigroup and $\Phi$ an isomorphism of $\cPSA(S)$ onto $\cPSA(T)$, we say that $\Phi$ is a {\em $\cPSA$-isomorphism} of $S$ onto $T$. It is clear that if $S$ is an inverse semigroup, then $\cPA(S)$ is an inverse submonoid of $\cPSA(S)$.
\br\label{33} {\rm (\cite[Lemma 2.3]{key3})} If $S$ and $T$ are inverse semigroups and $\Phi$ is a $\cPSA$-isomorphism of 
$S$ onto $T$, then $\Phi\vert_{\cPA(S)}$ is a $\cPA$-isomorphism of $S$ onto $T$. 
\er
Combining Theorem \ref{32} with a corollary of Theorem 4.13 of \cite{key5}, we can prove our final new result: 
\bt\label{34} Let $S$ be a tightly connected fundamental inverse semigroup with no nontrivial isolated subgroups, and let $T$ be an arbitrary semigroup. Then $\cPSA(S)\cong\cPSA(T)$ if and only if either $S\cong T$ or $(S, \leq)$ and $(T, \leq)$ are dually isomorphic chains.
\et 
{\bf Proof.} Suppose that $\cPSA(S)\cong\cPSA(T)$, and let $\Phi$ be a $\cPSA$-isomorphism of $S$ onto $T$. Since $S$ has no nontrivial isolated subgroups, by Corollary 4.14(b) of \cite{key5}, $T$ is an inverse semigroup with no nontrivial isolated subgroups. By Result \ref{33}, $\Phi\vert_{\cPA(S)}$ is a $\cPA$-isomorphism of $S$ onto $T$. Therefore, by Theorem \ref{32}, either $S\cong T$ or $(S, \leq)$ and $(T, \leq)$ are dually isomorphic chains. The converse is obvious.\epr\\
\sm
\section{Concluding remarks and open questions}
\sm Following Jones \cite{key10}, we will say that an inverse semigroup $S$ is {\em faintly archimedean} if whenever an idempotent $e$ of $S$ is strictly below every idempotent of a bicyclic or free inverse subsemigroup $\langle x\rangle$ of $S$, then $e < x$, and {\em quasi-archimedean} if it is faintly archimedean and $\langle x\rangle$ is combinatorial for every $x\in\Ng S$. (Thus a combinatorial inverse semigroup is quasi-archimedean if and only if it is faintly archimedean.) Jones proved \cite[Theorem 4.3]{key10} that if a combinatorial inverse semigroup $S$ is quasi-archimedean and $\Psi$ is a lattice isomorphism of $S$ onto an inverse semigroup $T$ such that the $E$-bijection $\psi_E$ is an isomorphism of $E_S$ onto $E_T$, then the base bijection $\psi$ is the unique isomorphism of $S$ onto $T$ inducing $\Psi$. Comparing this theorem with Result \ref{22} (which is Theorem 5 of \cite{key4}), it is natural to wonder how they are related. It was shown in \cite{key5} that for a combinatorial inverse semigroup $S$ neither of the two properties, shortly connected and quasi-archimedean, implies the other, so neither of the two theorems, Theorem 5 of \cite{key4} and Theorem 4.3 of \cite{key10}, is a corollary of the other one. 

It was shown in \cite{key10} (Theorem 4.5) that if $S$ is a completely semisimple, quasi-archimedean inverse semigroup, in which every noncombinatorial $\cD$-class contains at least {\em three} idempotents, and $\Psi$ is a lattice isomorphism of $S$ onto an inverse semigroup $T$ such that $\psi_E$ is an isomorphism of $E_S$ onto $E_T$, then the base bijection $\hat\psi$ is the unique isomorphism of $S$ onto $T$ inducing $\Psi$ (we express this briefly by saying that 
$S$ is {\em strongly lattice determined ``modulo semilattices''}). Now let $S$ be the semigroup constructed in the example preceding Lemma \ref{21}. Being finite, $S$ is quasi-archimedean and completely semisimple, but, since $D_{e_0}$ is its {\em noncombinatorial} $\cD$-class with just {\em two} idempotents, neither Theorem 4.3 nor Theorem 4.5 of \cite{key10} can be used for deciding whether $S$ is lattice determined. On the other hand, $S$ is a tightly connected fundamental inverse semigroup with no nontrivial isolated subgroups, so it is strongly lattice determined ``modulo semilattices'' by Theorem \ref{24} of this paper. (Of course, lattice determinability of this semigroup $S$ follows also from the theorem stated in Remark 2 on page 408 of \cite{key4}: {\em if $E$ is a semilattice such that the group of automorphisms of each principal ideal of $E$ is finite, then the Munn semigroup $T_E$ of $E$ is lattice determined.})

We conclude with two open questions. In view of Theorem \ref{24}, we would like to specialize Question 43.7 of 
\cite{key21} as follows: 
\bq\label{42}
Under what conditions is it true that, in the notation and under the assumptions of Theorem $\ref{24}$, the base bijection
$\hat\psi$, which is an isomorphism of $S$ onto $T$, induces $\Psi$ (and thus is the only isomorphism of $S$ onto $T$ inducing $\Psi$)?\eq 
Of course, finding an answer to Question \ref{42} will help answering a similar question regarding Theorem \ref{32}.
Finally, of particular interest to us is the following:\bq\label{43} Is it true that a fundamental inverse semigroup $S$ with no nontrivial isolated subgroups, which is shortly connected but not necessarily tightly connected, is lattice determined ``modulo semilattices'' (that is, under the assumption that the corresponding $E$-bijection is an isomorphism), and is it true that such a semigroup $S$ (assuming it is not a chain with respect to the natural order relation) is $\cPA$-determined?
\eq
\noindent {\bf Remark}. After this paper had been submitted for publication, Peter Jones noted that an inverse semigroup $S$ is tightly connected if and only if it is shortly connected and $\Ng S\cup E_S$ is its order ideal, giving an alternative characterization of tightly connected inverse semigroups.  
\med
\begin{center}
{\bf Acknowledgement}
\end{center}
\sm This article was written while I held the position of a Visiting Professor at the Department of Mathematics at
Stanford University from January through June of 2004. I am grateful to Professor Yakov Eliashberg for the invitation to Stanford and for his hospitality as well as to the department's faculty and staff for creating excellent working conditions. I am also grateful to the referee for useful comments and suggestions concerning the exposition of this paper.\\
%\sm

\end{document}